\documentclass{amsart}
\usepackage{amsfonts}
\usepackage{graphicx}
\usepackage{amscd}
\usepackage{amsmath}

\newtheorem{theorem}{Theorem}
\theoremstyle{plain}
\newtheorem{acknowledgement}{Acknowledgement}

\newtheorem{corollary}{Corollary}

\newtheorem{definition}{Definition}

\newtheorem{remark}{Remark}

\numberwithin{equation}{section}

\input{tcilatex}

\begin{document}
\title[$q$-Hardy-Berndt type sums ...]{$q$-Hardy-Berndt type sums associated
with $q$-Genocchi type zeta and $l$-functions}
\author{YILMAZ SIMSEK}
\address{AKDENIZ UNIVERSITY FACULTY OF ART\ AND\ SCIENCE DEPARTMENT OF
MATHEMATICS 07058\\
ANTALYA \\
TURKEY\\
Tel: ++90 242 310 23 43\\
Fax:+ +90242227 89 11}
\email{yilmazsimsek@hotmail.com, simsekyil63@yahoo.com}

\begin{abstract}
The aim of this paper is to define new generating functions. By applying the
Mellin transformation formula to these generating functions, we define $q$%
-analogue of Genocchi zeta function, $q$-analogue Hurwitz type Genocchi zeta
function, $q$-analogue Genocchi type $l$-function and two-variable $q$%
-Genocchi type $l$-function. Furthermore, we construct new genereting
functions of $q$-Hardy-Berndt type sums and $q$-Hardy-Berndt type sums
attached to Dirichlet character. We also give some new relations related to $%
q$-Hardy-Berndt type sums and $q$-Genocchi zeta function as well.
\end{abstract}

\maketitle

\section{Introduction}

In \cite{SimsekJMAA}, we defined generating functions. By using these
functions, we constructed $q$-Riemann zeta function, $q$-$L$-function and $q$%
-Dedekind type sum. This sum is defined by means of the generating function, 
$Y_{p}(h,k;q)$:%
\begin{equation*}
Y_{p}(h,k,q)=\sum_{m=1}^{\infty }\frac{f(-\frac{2mi\pi h}{k},q)-f(\frac{%
2mi\pi h}{k},q)}{m^{p}}
\end{equation*}%
where $h$ and $k$ are coprime positive integers and $p$ is an odd integer $%
\geq 1$, where 
\begin{equation}
f(t,q)=\sum_{n=1}^{\infty }q^{-n}\exp \left( -q^{-n}[n]t\right) \text{, cf. (%
\cite{Simsek3}, \cite{SimsekJMAA}).}  \label{f(t,q)}
\end{equation}

In the remainder of our work, we use\ $\exp (x)=e^{x}$ and $\chi $ is a
Dirichlet character of conductor $f\in \mathbb{Z}^{+}$, the set of positive
integer numbers.

The $q$-Dedekind type sum is given by\cite{SimsekJMAA}:

\begin{theorem}
\label{Theotem 1}Let $h$ and $k$ be positive integers and $(h,k)=1$ and
assume that $p$ is an odd integer $\geq 1$. We have%
\begin{equation*}
S_{p}(h,k;q)=\frac{p!}{\left( 2\pi i\right) ^{p}}Y_{p}(h,k;q).
\end{equation*}
\end{theorem}

By using (\ref{f(t,q)}), we construct $q$-analogous of Hardy-Berndt sums $%
s_{1}(h,k)$ and $s_{4}(h,k)$. Our aim is to define generating functions of\ $%
q$-analogous of the Hardy-Berndt type sums $S(h,k)$, $s_{2}(h,k)$, $%
s_{3}(h,k)$ and $s_{5}(h,k)$. We define%
\begin{equation}
F(t,q)=\sum_{n=1}^{\infty }(-1)^{n}q^{-n}\exp \left( -q^{-n}[n]t\right) .
\label{b3-1}
\end{equation}%
Applying the Mellin transformations to the equation (\ref{b3-1}), we define,
in Section 3, $q$-Genocchi type zeta functions, which implies the classical
Genocchi zeta functions. By applying the Mellin transformation to%
\begin{equation*}
F(t,x,q)=F(t,q)\exp (-tx),
\end{equation*}
we define Hurwitz type $q$-Genocchi zeta functions.

Character analogues of (\ref{b3-1}) is defined by: 
\begin{equation}
F_{\chi }(t,q)=\sum_{n=1}^{\infty }(-1)^{n}\chi (n)q^{-n}\exp \left(
-q^{-n}[n]t\right) .  \label{b3-5}
\end{equation}%
By applying the Mellin transformations to the equation (\ref{b3-5}), we
define $q$-Genocchi type $l$-function. By applying the Mellin transformation
to%
\begin{equation*}
F_{\chi }(t,x,q)=F_{\chi }(t,q)\exp (-tx),
\end{equation*}%
we define two-variable $q$-Genocchi type $l$-function. In Section 3, we
prove the relations between $q$-Genocchi type $l$-function, two-variable $q$%
-Genocchi type $l$-function\ and Hurwitz type $q$-Genocchi zeta functions.

By using (\ref{f(t,q)}), (\ref{b3-1}) and (\ref{b3-5}), we prove our main
results in Section 4-5. By using (\ref{b3-1}), $q$analogous of the
Hardy-Berndt type sums' $S(h,k;q)$ is defined by means of the generating
function, $Y_{0}(h,k;q)$:%
\begin{equation}
Y_{0}(h,k;q)=\sum_{m=1}^{\infty }\frac{F(-\frac{(2m-1)i\pi h}{2k},q)-F(\frac{%
(2m-1)i\pi h}{2k},q)}{2m-1}  \label{b4-1}
\end{equation}%
where $h$ and $k$ are coprime positive integers with $k\geq 1$.

\begin{theorem}
\label{Theorem-2}Let $h$ and $k$ be denote relatively prime integers with $%
k\geq 1$. If $h+k$ is odd, then%
\begin{equation}
S(h,k;q)=\frac{4}{\pi i}Y_{0}(h,k;q).  \label{A5}
\end{equation}
\end{theorem}

Theorem 2 implies the classical Hardy-Berndt sums $S(h,k)$. The other
theorems, which are related to $q$-Hardy-Bernd type sums, $s_{j}(h,k;q)$, $%
j=1,2,3,4,5$, are given in Section 4-5.

We define the sum $Y_{0,\chi }(t,q)$ as follows:%
\begin{equation}
Y_{0,\chi }(h,k;q)=\sum_{m=1}^{\infty }\frac{F_{\chi }(-\frac{(2m-1)i\pi h}{%
2k},q)-F_{\chi }(\frac{(2m-1)i\pi h}{2k},q)}{2m-1},  \label{A8}
\end{equation}%
where $h$ and $k$ are coprime positive integers. Therefore, generalized $q$%
-Hardy-Berndt type sums $S_{\chi }(h,k;q)$, with attached to $\chi $, are
given by the following theorem:

\begin{theorem}
\label{Thorem3}Let $h$ and $k$ be denote relatively prime integers with $%
k\geq 1$. If $h+k$ is odd, then%
\begin{equation}
S_{\chi }(h,k;q)=\frac{4}{\pi i}Y_{0,\chi }(h,k;q).  \label{bh0}
\end{equation}
\end{theorem}

\section{Definition and Notations}

The classical Dedekind sums $s(h,k)$ first arose in the transformation
formulae of the logarithm of the Dedekind eta-function. Similarly, the
Hardy-Berndt sums arose in the transformation formulae of the logarithm of
the theta-functions, $Log\vartheta _{n}(0,q)$, $n=2,3,4$\ cf. (\cite%
{berndt.Hardy}, \cite{Berndt and Goldberg1}, \cite{Goldberg}, \cite{Hardy}, 
\cite{Simsek2}).

Due to Hardy\cite{Hardy} and Berndt\cite{berndt.Hardy}, the Hardy-Berndt
sums are defined by:%
\begin{eqnarray*}
S(h,k) &=&\sum\limits_{j=1}^{k-1}(-1)^{j+1+[\frac{hj}{k}]},\text{ } \\
\text{\ \ \ \ \ \ \ \ \ \ \ \ \ \ \ \ }s_{1}(h,k)
&=&\sum\limits_{j=1}^{k}(-1)^{[\frac{hj}{k}]}((\frac{j}{k})), \\
s_{2}(h,k) &=&\sum\limits_{j=1}^{k}(-1)^{j}((\frac{j}{k}))((\frac{hj}{k})),%
\text{ } \\
\text{\ \ \ \ \ \ \ }s_{3}(h,k) &=&\sum\limits_{j=1}^{k}(-1)^{j}((\frac{hj}{k%
})), \\
s_{4}(h,k) &=&\sum\limits_{j=1}^{k-1}(-1)^{[\frac{hj}{k}]},\text{ \ \ } \\
\text{\ \ \ \ \ \ \ \ \ \ \ \ \ \ \ \ \ \ \ }s_{5}(h,k)
&=&\sum\limits_{j=1}^{k}(-1)^{j+[\frac{hj}{k}]}((\frac{j}{k})).
\end{eqnarray*}

Dieter\cite{Dieter} defined Hardy-Berndt sums by means of cotangent
function. Goldberg\cite{Goldberg} discovered some three-term and mixed
three-term relations for Hardy-Berndt sums. His proofs are based on Berndt's
transformation formulae for the logarithms of the classical theta-functions.
For an elaboration of this connection and fundamental properties of
Hardy-Berndt sums cf. (see for detail \cite{Apostol and Vu}, \cite{Berndt2}, 
\cite{Berndt and Dieter1}, \cite{Dieter}, \cite{Petet and Sitara}, \cite%
{simsekTurJM}, \cite{Yilmaz1}).

Hardy-Berndt sums were represented as infinite trigonometric sums by Berndt
and Goldberg\cite{Berndt and Goldberg1} as follows:

\begin{theorem}
\label{Theorem 4}Let $h$ and $k$ denote relatively prime integers with $k>0$%
. If $h+k$ is odd, then 
\begin{equation}
S(h,k)=\frac{4}{\pi }\sum_{n=1}^{\infty }\frac{\tan (\frac{\pi h(2n-1)}{2k})%
}{2n-1},  \label{Eq-1}
\end{equation}

if $h$ is even and $k$ is odd, then 
\begin{equation}
s_{1}(h,k)=-\frac{2}{\pi }\sum_{%
\begin{array}{c}
n=1 \\ 
2n-1\not\equiv 0(\func{mod}k)%
\end{array}
}^{\infty }\frac{\cot (\frac{\pi h(2n-1)}{2k})}{2n-1},  \label{Eq-2}
\end{equation}

if $h$ is odd and $k$ is even, then 
\begin{equation}
s_{2}(h,k)=-\frac{1}{2\pi }\sum_{%
\begin{array}{c}
n=1 \\ 
2n\not\equiv 0(\func{mod}k)%
\end{array}%
}^{\infty }\frac{\tan (\frac{\pi hn)}{k})}{n},  \label{Eq-3}
\end{equation}%
if $k$ is odd, then 
\begin{equation}
s_{3}(h,k)=\frac{1}{\pi }\sum_{n=1}^{\infty }\frac{\tan (\frac{\pi hn}{k})}{n%
},  \label{Eq-4}
\end{equation}%
if $h$ is odd, then 
\begin{equation}
s_{4}(h,k)=\frac{4}{\pi }\sum_{n=1}^{\infty }\frac{\cot (\frac{\pi h(2n-1)}{%
2k})}{2n-1},  \label{Eq-5}
\end{equation}%
and if $h$ and $k$ are odd, then 
\begin{equation}
s_{5}(h,k)=\frac{2}{\pi }\sum_{%
\begin{array}{c}
n=1 \\ 
2n-1\not\equiv 0(\func{mod}k)%
\end{array}%
}^{\infty }\frac{\tan (\frac{\pi h(2n-1)}{2k})}{2n-1}.  \label{Eq-6}
\end{equation}
\end{theorem}

Observe that if $q\rightarrow 1$, then Theorem 2 reduces to (\ref{Eq-1}),
which are stating in Section 4 in detail.

In \cite{Sitaramachandrarao}, Sitaramachandrarao studied on Hardy-Berndt
sums. He proved Reciprocity Law of these sums and gave the relations between
Hardy-Berndt sums and Dedekind sums. In \cite{Simsek.JKMS}, we defined $p$%
-adic Hardy-type sums and $p$-adic $q$-higher-order Hardy-type sums. We gave 
$p$-adic continuous functions. By using these functions and $p$-adic $q$%
-integral, we obtain relations between $p$-adic $q$-higher-order Hardy-type
sums, Bernoulli functions and Lambert series.

If $q\in \mathbb{C}$, the field of complex numbers, then we assume $\mid
q\mid <1$. We set 
\begin{equation*}
\lbrack x]=[x:q]=\frac{1-q^{x}}{1-q}.
\end{equation*}%
\ Note that $\lim_{q\rightarrow 1}[x]=x$, cf. (\cite{TKim-4}, \cite{Simsek3}%
, \cite{Simsek4}).

The Euler numbers $E_{n}$ are usually defined by means of of the following
generating function cf. (\cite{kim-kitap}, \cite{Kim-genochi}, \cite%
{Shiratani}, \cite{Srivastava-Kim-Sim}): 
\begin{equation*}
\frac{2}{\exp (t)+1}=\sum_{n=0}^{\infty }E_{n}\frac{t^{n}}{n!}\text{, }%
|t|<\pi .
\end{equation*}%
The Genocchi numbers $G_{n}$ are usually defined by means of of the
following generating function cf. (\cite{kim-kitap}, \cite{Kim-genochi}, 
\cite{kimsimnaci}): 
\begin{equation*}
\frac{2t}{\exp (t)+1}=\sum_{n=0}^{\infty }G_{n}\frac{t^{n}}{n!}\text{, }%
|t|<\pi .
\end{equation*}%
These numbers are classical and important in number theory. In \cite%
{Kim-genochi}, Kim defined generating function of the $q$-Genocchi numbers
and $q$-Euler numbers as follows%
\begin{equation*}
\lbrack 2]\exp (\frac{t}{1-q})\sum_{n=0}^{\infty }\frac{(-1)^{n}}{%
(1+q^{n+1})(1-q)^{n}}\frac{t^{n}}{n!}=\sum_{m=0}^{\infty }E_{m,q}\frac{t^{m}%
}{m!},
\end{equation*}%
where $E_{m,q}$ is denoted $q$-Euler numbers.%
\begin{eqnarray*}
G_{q}(t) &=&[2]t\sum_{m=0}^{\infty }(-1)^{n}q^{n}e^{[n]t} \\
&=&\sum_{m=0}^{\infty }G_{m,q}\frac{t^{m}}{m!},
\end{eqnarray*}%
where $G_{m,q}$ is denoted $q$-Genocchi numbers.

Genocchi zeta function is defined by cf. (\cite{kim-kitap}, p. 108):%
\begin{equation*}
\zeta _{G}(s)=2\sum_{n=1}^{\infty }\frac{(-1)^{n}}{n^{s}},
\end{equation*}%
where $s\in \mathbb{C}$.

In \cite{kimsimnaci}, by using q-Volkenborn Integral, Kim defined generating
functions. By applying this generating function, they constructed $q$%
-Genocchi zeta function, the others $q$-function. They gave relations
between $q$-Genocchi zeta function and $q$-Genocchi numbers. They also
defined high order of $q$-Genocchi zeta function and $q$-Genocchi numbers
and polynomials.

Kim and Rim\cite{Tkim-SHRim} defined two-variable $L$-function. They gave
main properties of this function. In \cite{TKim-4}, Kim constructed the
two-variable $p$-adic $q$-$L$-function which interpolates the generalized $q$%
-Bernoulli polynomials attached to Dirichlet character. In \cite%
{Y.Simsek-Dkim-SH. Rim}, the author, Kim and Rim constructed the
two-variable Dirichlet $q$-$L$-function and the two-variable multiple
Dirichlet-type Changhee $q$-$L$-function.

We summarize our paper as follows:

In Section 3, we give new generating functions. By applying the Mellin
transformation formula to these generating functions, we will define $q$%
-analogue of Genocchi zeta function, $q$-analogue Hurwitz type Genocchi zeta
function, $q$-analogue $l$-function and two-variable $q$-$l$-function. In
Section 4, by using generating functions in Section 3, we will construct new
generating function which produce $q$-Hardy-Berndt type sums. In Section 5,
by using generating functions in Section 3, we will construct new generating
function attached to Dirichlet character which produces $q$-Hardy-Berndt
type sums attached to Dirichlet character. In section 6, we define new
generating functions. By applying Mellin transformation to these functions,
we can give some new relations which are related to Riemann zeta functions, $%
q$-Riemann zeta function, $q$-$L$-function and $q$-Genocchi zeta function.

\section{$q$-Genocchi zeta function and $l$-function}

In \cite{Simsek1}, \cite{Simsek3}, and \cite{Simsek4}, the author defined
generating functions, which are interpolates twisted Bernoulli numbers and
polynomials, twisted Euler numbers and polynomials. In this section, we give
new generating functions which produce $q$-Genocchi zeta functions and $q$-$%
l $-series with attached to Dirichlet character. Therefore, by using these
generating functions, we construct new $q$-analogue of Hardy-Berndt sums in
the next sections. We also give relations between these sums, $q$-Genocchi
zeta functions and $q$-$l$-series.

\begin{remark}
\begin{equation*}
\zeta _{G}(s)=\sum_{n=1}^{\infty }\frac{(-1)^{n}}{n^{s}}
\end{equation*}%
and%
\begin{equation}
\zeta _{G}(s)\Gamma (s)=\int_{0}^{\infty }\frac{2x^{s-1}}{\exp (-x)+1}dx,
\label{Eq-12}
\end{equation}%
where $\Gamma (s)$ is Euler's gamma function\ and 
\begin{equation*}
\zeta _{G}(1-n)=-\frac{G_{n}}{n}\text{, \ }n>1\text{ cf. (\cite%
{jangkimLeePark}, \cite{kim-kitap}, p. 108, Eq. (2.43)).}
\end{equation*}

In (\ref{Eq-12}), due to \cite{kim-kitap}, \cite{Waldschmidt} and \cite%
{Simsek-Yang}, 
\begin{equation*}
\sum_{n=1}^{\infty }(-1)^{n}\int_{0}^{\infty }x^{s-1}\exp (-nx)dx,
\end{equation*}%
and if $q\rightarrow 1$ in (\ref{b3-1}), we have 
\begin{align*}
\lim_{q\rightarrow 1}2F(t,q)& =2\sum_{n=1}^{\infty }(-1)^{n}\exp (-nt) \\
& =\frac{2}{\exp (t)+1}.
\end{align*}%
According to the definition of Genocchi numbers: 
\begin{equation*}
\frac{2}{\exp (t)+1}=\sum_{n=0}^{\infty }G_{n}\frac{t^{n}}{n!}\text{, }%
|t|<\pi
\end{equation*}%
we get an asymptotic expansion near 0 
\begin{equation*}
\frac{1}{\exp (t)+1}\sim \sum_{k=-1}^{\infty }G_{k+1}\frac{t^{k}}{(k+1)!}
\end{equation*}%
while 
\begin{equation*}
\frac{1}{\exp (t)+1}\sim 0
\end{equation*}%
for $t$ near $\infty $.
\end{remark}

Applying the Mellin transformations to (\ref{b3-1}), we find that%
\begin{align*}
& \frac{1}{\Gamma (s)}\int_{0}^{\infty }t^{s-1}F(t,q)dt \\
& =\frac{1}{\Gamma (s)}\int_{0}^{\infty }t^{s-1}\left( \sum_{n=1}^{\infty
}(-1)^{n}q^{-n}\exp \left( -q^{-n}[n]t\right) \right) dt \\
& =\frac{1}{\Gamma (s)}\sum_{n=1}^{\infty }(-1)^{n}q^{-n}\int_{0}^{\infty
}t^{s-1}\exp \left( -q^{-n}[n]t\right) dt \\
& =\frac{1}{\Gamma (s)}\sum_{n=1}^{\infty }\frac{(-1)^{n}q^{-n}}{\left(
q^{-n}[n]\right) ^{s}}\int_{0}^{\infty }u^{s-1}e^{-u}du \\
& =\sum_{n=1}^{\infty }\frac{(-1)^{n}q^{-n}}{\left( q^{-n}[n]\right) ^{s}}%
=\Im _{q}(s).
\end{align*}%
The right-hand side of the above converges when $\func{Re}(s)>1$. By using
the above series, we are now ready to define $q$-analogue of the Genocchi
zeta functions.

\begin{definition}
Let $s\in \mathbb{C}$ and $\func{Re}(s)>1$. $q$-analogue of the Genocchi
type zeta function expressed by the formula 
\begin{equation}
\Im _{G,q}(s)=[2]\Im _{q}(s).  \label{b3-2}
\end{equation}
\end{definition}

\begin{remark}
Observe that when $q\rightarrow 1$, (\ref{b3-2}) reduces to ordinary
Genocchi zeta functions (see \cite{kim-kitap}, p. 108). In \cite%
{cenkci-can-kurt}, Cenkci, Can and Kurt defined different type $q$-Genocch
zeta functions, which is defined as follows%
\begin{equation*}
\zeta _{q}^{(G)}(s)=q(1+q)\sum_{n=1}^{\infty }\frac{(-1)^{n+1}q^{n}}{[n]^{s}}%
.
\end{equation*}
\end{remark}

We define $q$-analogue of the Hurwitz type Genocchi zeta function by means
of the generating function%
\begin{equation}
F(t,x,q)=F(t,q)\exp (-tx)=\sum_{n=0}^{\infty }(-1)^{n}q^{-n}\exp \left(
-\left( q^{-n}[n]+x\right) t\right) .  \label{b3-3}
\end{equation}

By applying the Mellin transformations to (\ref{b3-3}), we obtain%
\begin{align*}
& \frac{1}{\Gamma (s)}\int_{0}^{\infty }t^{s-1}F(t,x,q)dt \\
& =\frac{1}{\Gamma (s)}\int_{0}^{\infty }t^{s-1}\left( \sum_{n=0}^{\infty
}(-1)^{n}q^{-n}\exp \left( -\left( q^{-n}[n]+x\right) t\right) \right) dt \\
& =\frac{1}{\Gamma (s)}\sum_{n=0}^{\infty }(-1)^{n}q^{-n}\int_{0}^{\infty
}t^{s-1}\exp \left( -\left( q^{-n}[n]+x\right) t\right) dt \\
& =\frac{1}{\Gamma (s)}\sum_{n=0}^{\infty }\frac{(-1)^{n}q^{-n}}{\left(
q^{-n}[n]+x\right) ^{s}}\int_{0}^{\infty }u^{s-1}e^{-u}du \\
& =\sum_{n=0}^{\infty }\frac{(-1)^{n}q^{-n}}{\left( q^{-n}[n]+x\right) ^{s}}%
=\Im _{q}(s,x).
\end{align*}%
By using the above equation, we are ready to define $q$-analogue of the
Hurwitz-type Genocchi zeta functions.

\begin{definition}
Let $s\in \mathbb{C}$ and $\func{Re}(s)>1$ and $0<x\leq 1$. $q$-analogue of
the Hurwitz-type Genocchi zeta function expressed by the formula 
\begin{equation}
\Im _{G,q}(s,x):=[2]\Im _{q}(s,x).  \label{b3-4}
\end{equation}
\end{definition}

\begin{remark}
We give another version of (\ref{b3-3}) as follows:%
\begin{eqnarray*}
F(t,x,q) &=&F(t,q)\exp (-t[x])=\sum_{n=0}^{\infty }(-1)^{n}q^{-n}\exp \left(
-\left( q^{-n}[n]+[x]\right) t\right) \\
&=&\sum_{n=0}^{\infty }(-1)^{n}q^{-n}\exp \left( -q^{-n}t([n]+q^{n}[x]\right)
\end{eqnarray*}%
By using well-known the identity%
\begin{equation*}
\lbrack n+x]=[n]+q^{n}[x]
\end{equation*}%
we have%
\begin{equation*}
F(t,x,q)=\sum_{n=0}^{\infty }(-1)^{n}q^{-n}\exp \left( -q^{-n}t[n+x]\right) .
\end{equation*}%
By applying the Mellin transformation to the above generating function, we
obtain%
\begin{eqnarray*}
\frac{1}{\Gamma (s)}\int_{0}^{\infty }t^{s-1}F(t,x,q)dt
&=&\sum_{n=0}^{\infty }\frac{(-1)^{n}q^{-n}}{\left( q^{-n}[n]+[x]\right) ^{s}%
} \\
&=&\sum_{n=0}^{\infty }\frac{(-1)^{n}q^{-n(1-s)}}{[n+x]^{s}}=\Im _{q}(s,x)
\end{eqnarray*}%
Observe that if $x=1$, then $\Im _{G,q}(s,x)$\ is reduced $\Im _{G,q}(s)$
and if $q\rightarrow 1$, then $\Im _{G,q}(s,x)\rightarrow \Im _{G}(s,x)$. A
function $\Im _{G}(s,x)$ is called a ordinary Hurwitz-type Genocchi zeta
function if $\Im _{G}(s,x)$ is expressed by the formula%
\begin{equation*}
\Im (s,x):=2\sum_{n=0}^{\infty }\frac{(-1)^{n}}{\left( n+x\right) ^{s}},
\end{equation*}%
where $s\in \mathbb{C}$, $\func{Re}(s)>1$ and $0<x\leq 1$ (see \cite%
{kimsimnaci}).
\end{remark}

By applying the Mellin transformations to the equation (\ref{b3-5}), we
obtain%
\begin{align*}
& \frac{1}{\Gamma (s)}\int_{0}^{\infty }t^{s-1}F_{\chi }(t,q)dt \\
& =\frac{1}{\Gamma (s)}\int_{0}^{\infty }t^{s-1}\left( \sum_{n=1}^{\infty
}(-1)^{n}\chi (n)q^{-n}\exp \left( -q^{-n}[n]t\right) \right) dt \\
& =\frac{1}{\Gamma (s)}\sum_{n=1}^{\infty }(-1)^{n}\chi
(n)q^{-n}\int_{0}^{\infty }t^{s-1}\exp \left( -q^{-n}[n]t\right) dt \\
& =\frac{1}{\Gamma (s)}\sum_{n=1}^{\infty }\frac{(-1)^{n}\chi (n)q^{-n}}{%
\left( q^{-n}[n]\right) ^{s}}\int_{0}^{\infty }u^{s-1}e^{-u}du \\
& =\sum_{n=1}^{\infty }\frac{(-1)^{n}\chi (n)q^{-n}}{\left( q^{-n}[n]\right)
^{s}}=l_{q}(s,\chi ).
\end{align*}%
Now, by using the above equation, we define $q$-analogue (Genocchi-type) $l$%
-function as follows:

\begin{definition}
Let $\chi $ be a Dirichlet character. Let $s\in \mathbb{C}$ and $\func{Re}%
(s)>1$. (Genocchi-type) $q$-$l$-function expressed by the formula 
\begin{equation}
l_{G,q}(s,\chi )=[2]l_{q}(s,\chi ).  \label{b3-6}
\end{equation}
\end{definition}

A function $l_{G}(s,\chi )$ is called a ordinary Genocchi-type $l$-function
if $l_{G}(s,\chi )$ is expressed by the formula%
\begin{equation*}
l(s,x):=2\sum_{n=0}^{\infty }\frac{(-1)^{n}\chi (n)}{\left( n+x\right) ^{s}},
\end{equation*}%
where $s\in \mathbb{C}$, $\func{Re}(s)>1$ and $0<x\leq 1$ \cite{kimsimnaci}.

Observe that when $\chi \equiv 1$, (\ref{b3-6}) reduces to (\ref{b3-2}):%
\begin{equation*}
l_{q}(s,1)=\Im _{q}(s).
\end{equation*}%
In \cite{cenkci-can-kurt}, Cenkci, Can and Kurt defined Genocch measure. By
using this measure and Volkenborn Integral, they defined different type $l$%
-function. This functions interpolate generalized Genocchi numbers.

Relation between $\Im _{G,q}(s,x)$ and $l_{G,q}(s,\chi )$ and is given as
follows:

\begin{theorem}
\label{Theorem 5}Let $\chi $ be a Dirichlet character. We have%
\begin{equation}
l_{G,q}(s,\chi )=\frac{1}{[f]^{s}}\sum_{a=1}^{f}(-1)^{a}q^{a(s-1)}\chi
(a)\Im _{G,q^{f}}\left( s,\frac{[a]}{[f]}\right) .  \label{b3-7}
\end{equation}
\end{theorem}

\begin{proof}
Substituting $n=a+mf$, where $m=0,1,2,3,...,\infty $ and $a=1,...,f$ into (%
\ref{b3-6}), we obtain 
\begin{equation*}
l_{G,q}(s,\chi )=\sum_{a=1}^{f}(-1)^{a}q^{-a}\chi (a)\sum_{m=0}^{\infty }%
\frac{(-1)^{mf}q^{-mf}}{(q^{-a-mf}[a+mf])^{s}}.
\end{equation*}%
By using $[a+mf]=[m:q^{f}][f]+q^{mf}[a]$ \ in the above equation, we obtain%
\begin{eqnarray*}
l_{G,q}(s,\chi ) &=&\sum_{a=1}^{f}(-1)^{a}q^{-a}\chi (a)\sum_{m=0}^{\infty }%
\frac{(-1)^{mf}q^{-mf}}{\left( q^{-a-mf}\left( [m:q^{f}][f]+q^{mf}[a]\right)
\right) ^{s}} \\
&=&\frac{1}{[f]^{s}}\sum_{a=1}^{f}(-1)^{a}q^{a(s-1)}\chi
(a)\sum_{m=0}^{\infty }\frac{(-1)^{mf}q^{-mf}}{\left( q^{-mf}[m:q^{f}]+\frac{%
[a]}{[f]}\right) ^{s}}
\end{eqnarray*}%
By using Eq.(\ref{b3-4}) in the above and after elementary calculations, we
easily arrive at (\ref{b3-7}).
\end{proof}

Now, we define the following generating function%
\begin{equation}
F_{\chi }(t,x,q)=F_{\chi }(t,q)e^{-tx}=\sum_{n=0}^{\infty }(-1)^{n}\chi
(n)q^{-n}\exp \left( -(q^{-n}[n]+x)t\right) .  \label{b3-8}
\end{equation}%
By using the Mellin transformations in the above equation, we obtain%
\begin{align*}
& \frac{1}{\Gamma (s)}\int_{0}^{\infty }t^{s-1}F_{\chi }(t,x,q)dt \\
& =\frac{1}{\Gamma (s)}\int_{0}^{\infty }t^{s-1}\left( \sum_{n=0}^{\infty
}(-1)^{n}\chi (n)q^{-n}\exp \left( -(q^{-n}[n]+x)t\right) \right) dt \\
& =\frac{1}{\Gamma (s)}\sum_{n=0}^{\infty }(-1)^{n}\chi
(n)q^{-n}\int_{0}^{\infty }t^{s-1}\exp \left( -(q^{-n}[n]+x)t\right) dt \\
& =\frac{1}{\Gamma (s)}\sum_{n=0}^{\infty }\frac{(-1)^{n}\chi (n)q^{-n}}{%
\left( q^{-n}[n]+x\right) ^{s}}\int_{0}^{\infty }u^{s-1}e^{-u}du \\
& =\sum_{n=0}^{\infty }\frac{(-1)^{n}\chi (n)q^{-n}}{\left(
q^{-n}[n]+x\right) ^{s}}=l_{q}(s,x,\chi ).
\end{align*}%
By using the above equation, we define the two-variable $q$-analogue $l$%
-function.

\begin{definition}
\begin{equation}
l_{G,q}(s,x,\chi )=[2]l_{q}(s,x,\chi ).  \label{b3-9}
\end{equation}
\end{definition}

Relation between $\Im _{G,q}(s,x)$ and $l_{G,q}(s,x,\chi )$ and is given as
follows

\begin{theorem}
\label{Theorem 6}Let $\chi $ be a Dirichlet character of conductor $f$. We
have%
\begin{equation}
l_{G,q}(s,x,\chi )=\frac{1}{[f]^{s}}\sum_{a=1}^{f}(-1)^{a}q^{a(s-1)}\chi
(a)\Im _{G,q^{f}}\left( s,\frac{[a]+xq^{a}}{[f]}\right) .  \label{aa.0}
\end{equation}
\end{theorem}

\begin{proof}
The proof of (\ref{aa.0}) follows precisely along the same lines as that of (%
\ref{b3-7}), and so we omit it.
\end{proof}

\section{$q$-Hardy-Berndt type sums}

In this section, we establish a general theorem related to the $q$%
-Hardy-Berndt type sums.

\begin{theorem}
\label{Theorem 7}Let $h$ and $k$ be coprime positive integers with $h\geq 1$%
. We have%
\begin{equation*}
Y_{0}(h,k;q)=2i\sum_{m=1}^{\infty }\sum_{n=1}^{\infty }\frac{%
(-1)^{n}q^{-n}\sin (\frac{q^{-n}[n](2m-1)\pi h}{2k})}{2m-1}.
\end{equation*}
\end{theorem}

\begin{proof}
By using (\ref{b3-1}) and (\ref{b4-1}), we have 
\begin{align*}
Y_{0}(h,k,q)& =\sum_{m=1}^{\infty }\frac{1}{2m-1}(\sum_{n=1}^{\infty
}(-1)^{n}q^{-n}\exp (\frac{q^{-n}[n](2m-1)\pi hi}{2k}) \\
& -\sum_{n=1}^{\infty }(-1)^{n}q^{-n}\exp (-\frac{q^{-n}[n](2m-1)\pi hi}{2k}%
)) \\
& =\sum_{m=1}^{\infty }\frac{1}{2m-1}\sum_{n=1}^{\infty }(-1)^{n}q^{-n}(\exp
(\frac{q^{-n}[n](2m-1)\pi hi}{2k}) \\
& -\exp (-\frac{q^{-n}[n](2m-1)\pi hi}{2k})).
\end{align*}%
Recalling that $2i\sin x=\exp (ix)-\exp (-ix)$, we easily complete the proof.
\end{proof}

We are now ready to prove the primary theorems of this section.

\begin{proof}[Proof of Theorem 2]
Applying the generating function (\ref{b4-1}), and using (\ref{b3-1}) and
Theorem 4, Eq. (\ref{Eq-1}) and recalling the definition of $S(h,k)$, after
some elementary calculations, we arrive at the desired result.
\end{proof}

Observe that when $q\rightarrow 1$ in Theorem 2, then we have%
\begin{eqnarray*}
Y_{0}(h,k;1) &=&\sum_{m=1}^{\infty }\frac{1}{2m-1}\sum_{n=1}^{\infty }(\exp (%
\frac{q^{-n}[n](2m-1)\pi hi}{2k}) \\
&&-\exp (-\frac{q^{-n}[n](2m-1)\pi hi}{2k})).
\end{eqnarray*}%
Hence, applying the well-known series%
\begin{equation*}
\sum_{m=1}^{\infty }(-1)^{n}x^{n}=\frac{x}{x+1}
\end{equation*}%
in the above equation, we deduce that%
\begin{eqnarray}
Y_{0}(h,k;1) &=&\sum_{m=1}^{\infty }\frac{1}{2m-1}(\frac{\exp (\frac{\pi
ih(2m-1)}{2k})}{1+\exp (\frac{\pi ih(2m-1)}{2k})}  \label{A6} \\
&&-\frac{\exp (-\frac{\pi ih(2m-1)}{2k})}{1+\exp (-\frac{\pi ih(2m-1)}{2k})}%
).  \notag
\end{eqnarray}%
We define%
\begin{equation}
i\tan (x\pi )=\frac{\exp (2ix\pi )}{1+\exp (2ix\pi )}-\frac{\exp (-2ix\pi )}{%
1+\exp (-2ix\pi )}.  \label{A5A}
\end{equation}%
By (\ref{A5}), (\ref{A6}) and (\ref{A5A}), we arrive at the following
corollary.

\begin{corollary}
Let $h$ and $k$ be denote relatively prime integers with $k\geq 1$. If $h+k$
is odd, then we have%
\begin{eqnarray*}
S(h,k;1) &=&\frac{4}{\pi i}Y_{0}(h,k;1) \\
&=&\frac{4}{\pi i}\sum_{m=1}^{\infty }\frac{1}{2m-1}(\frac{\exp (\frac{\pi
ih(2m-1)}{2k})}{1+\exp (\frac{\pi ih(2m-1)}{2k})} \\
&&-\frac{\exp (-\frac{\pi ih(2m-1)}{2k})}{1+\exp (-\frac{\pi ih(2m-1)}{2k})}%
).
\end{eqnarray*}
\end{corollary}

\begin{remark}
The case $q\rightarrow 1$, $S(h,k,1)$ is denoted the classical Hardy-Berndt
sums, which is given (\ref{Eq-1}). Generalized Dedekind sums, $s(h,k;p)$ are
expressible as infinite series related to certain Lambert series. A
representation of $s(h,k;p)$ as infinite series and trigonometric
representation was given by Apostol (\cite{Apostol1}, \cite{Apostol2}).\ $q$%
-Dedekind\ type sums given by the author\cite{SimsekJMAA}.
\end{remark}

We now define generating function of $q$-Hardy-Bernd type sums $s_{j}(h,k;q)$%
, $j=1,2,3,4,5$. The definition of $q$-analogous of $s_{j}(h,k)$, $%
j=1,2,3,4,5$ follow precisely along the same lines as the definition of the $%
q$-analogous of the Dedekind sums, which was established by the author\cite%
{SimsekJMAA}.

$q$-Hardy-Bernd type sums $s_{1}(h,k;q)$ is defined by means of the
generating function, $Y_{1}(h,k;q)$:%
\begin{equation}
Y_{1}(h,k;q)=\sum_{%
\begin{array}{c}
m=1 \\ 
2m-1\not\equiv 0\func{mod}k%
\end{array}%
}^{\infty }\frac{f(-\frac{(2m-1)i\pi h}{2k},q)-f(\frac{(2m-1)i\pi h}{2k},q)}{%
2m-1},  \label{S1(h,k)}
\end{equation}%
where $f(t,q)$ is defined by (\ref{f(t,q)}).

By using (\ref{S1(h,k)}), (\ref{f(t,q)}) and (\ref{Eq-2}), we obtain the
fallowing theorem.

\begin{theorem}
\label{Theorem 8}Let $h$ and $k$ be denote relatively prime integers with $%
k\geq 1$. If $h$ is even and $k$ is odd, then%
\begin{equation*}
s_{1}(h,k;q)=-\frac{2}{\pi i}Y_{1}(h,k;q).
\end{equation*}
\end{theorem}

\begin{corollary}
Let $h$ and $k$ be denote relatively prime integers with $k\geq 1$. If $h$
is even and $k$ is odd, then%
\begin{eqnarray*}
s_{1}(h,k;1) &=&-\frac{2}{\pi i}Y_{1}(h,k;1) \\
&=&-\frac{2}{\pi i}\sum_{%
\begin{array}{c}
m=1 \\ 
2m-1\not\equiv 0\func{mod}k%
\end{array}%
}^{\infty }\frac{1}{2m-1}(\frac{\exp (\frac{\pi ih(2m-1)}{2k})}{1-\exp (%
\frac{\pi ih(2m-1)}{2k})} \\
&&-\frac{\exp (-\frac{\pi ih(2m-1)}{2k})}{1-\exp (-\frac{\pi ih(2m-1)}{2k})}%
).
\end{eqnarray*}
\end{corollary}

Observe that the case $q\rightarrow 1$, $s_{1}(h,k;1)$ is the classical
Hardy-Berndt sums, which is given (\ref{Eq-2}).

By using (\ref{b3-1}), $q$-Hardy-Berndt type sums' $s_{2}(h,k;q)$ is defined
by means of the generating function, $Y_{2}(h,k;q)$:%
\begin{equation}
Y_{2}(h,k;q)=\sum_{%
\begin{array}{c}
m=1 \\ 
2m\not\equiv 0\func{mod}k%
\end{array}%
}^{\infty }\frac{F(-\frac{mi\pi h}{k},q)-F(\frac{mi\pi h}{k},q)}{m}
\label{s2(h,k)}
\end{equation}%
By using (\ref{s2(h,k)}), (\ref{b3-1}) and (\ref{Eq-3}), we obtain the
fallowing theorem.

\begin{theorem}
\label{Theorem 9}Let $h$ and $k$ be denote relatively prime integers with $%
k\geq 1$. If $h$ is odd and $k$ is even, then%
\begin{equation*}
s_{2}(h,k;q)=-\frac{1}{2\pi i}Y_{2}(h,k;q).
\end{equation*}
\end{theorem}

\begin{corollary}
Let $h$ and $k$ be denote relatively prime integers with $k\geq 1$. If $h$
is odd and $k$ is even, then%
\begin{eqnarray*}
s_{2}(h,k;1) &=&-\frac{1}{2\pi i}Y_{2}(h,k;1) \\
&=&-\frac{2}{\pi i}\sum_{%
\begin{array}{c}
m=1 \\ 
2m\not\equiv 0\func{mod}k%
\end{array}%
}^{\infty }\frac{1}{m}(\frac{\exp (\frac{\pi ihm}{k})}{1+\exp (\frac{\pi ihm%
}{k})} \\
&&-\frac{\exp (-\frac{\pi ihm}{k})}{1+\exp (-\frac{\pi ihm}{k})}).
\end{eqnarray*}
\end{corollary}

Note that if $q\rightarrow 1$, then $s_{2}(h,k;1)$ is denoted the classical
Hardy-Berndt sums, which is given (\ref{Eq-3}).

By using (\ref{b3-1}), $q$-Hardy-Berndt type sums' $s_{3}(h,k;q)$ is defined
by means of the generating function, $Y_{3}(h,k;q)$:%
\begin{equation}
Y_{3}(h,k;q)=\sum_{m=1}^{\infty }\frac{F(-\frac{mi\pi h}{k},q)-F(\frac{mi\pi
h}{k},q)}{m}  \label{s3(h,k)}
\end{equation}%
By using (\ref{s3(h,k)}), (\ref{b3-1}) and (\ref{Eq-4}), we obtain the
fallowing theorem.

\begin{theorem}
\label{Theorem 10}Let $h$ and $k$ be denote relatively prime integers with $%
k\geq 1$. If $k$ is odd, then%
\begin{equation*}
s_{3}(h,k;q)=\frac{1}{\pi i}Y_{3}(h,k;q).
\end{equation*}
\end{theorem}

\begin{corollary}
Let $h$ and $k$ be denote relatively prime integers with $k\geq 1$. If $k$
is odd, then 
\begin{eqnarray*}
s_{3}(h,k;1) &=&\frac{1}{\pi i}Y_{3}(h,k;1) \\
&=&\frac{1}{\pi i}\sum_{m=1}^{\infty }\frac{1}{m}(\frac{\exp (\frac{\pi ihm}{%
k})}{1+\exp (\frac{\pi ihm}{k})} \\
&&-\frac{\exp (-\frac{\pi ihm}{k})}{1+\exp (-\frac{\pi ihm}{k})}).
\end{eqnarray*}
\end{corollary}

Observe that the case $q\rightarrow 1$, $s_{3}(h,k;1)$ is denoted the
classical Hardy-Berndt sums, which is given (\ref{Eq-4}).

$q$-Hardy-Bernd type sums $s_{4}(h,k;q)$ is defined by means of the
generating function, $Y_{4}(h,k;q)$:%
\begin{equation}
Y_{4}(h,k;q)=\sum_{m=1}^{\infty }\frac{f(-\frac{(2m-1)i\pi h}{2k},q)-f(\frac{%
(2m-1)i\pi h}{2k},q)}{2m-1},  \label{s4(h,k)}
\end{equation}%
where $f(t,q)$ is defined by (\ref{f(t,q)}).

By using (\ref{s4(h,k)}), (\ref{f(t,q)}) and (\ref{Eq-5}), we obtain the
fallowing theorem.

\begin{theorem}
\label{Theorem 11}Let $h$ and $k$ be denote relatively prime integers with $%
k\geq 1$. If $h$ is odd, then%
\begin{equation*}
s_{4}(h,k;q)=\frac{4}{\pi i}Y_{4}(h,k;q).
\end{equation*}
\end{theorem}

\begin{corollary}
Let $h$ and $k$ be denote relatively prime integers with $k\geq 1$. If $h$
is odd, then%
\begin{eqnarray*}
s_{4}(h,k;1) &=&\frac{4}{\pi i}Y_{4}(h,k;1) \\
&=&\frac{4}{\pi i}\sum_{m=1}^{\infty }\frac{1}{2m-1}(\frac{\exp (\frac{\pi
ih(2m-1)}{2k})}{1-\exp (\frac{\pi ih(2m-1)}{2k})} \\
&&-\frac{\exp (-\frac{\pi ih(2m-1)}{2k})}{1-\exp (-\frac{\pi ih(2m-1)}{2k})}%
).
\end{eqnarray*}
\end{corollary}

Observe that the case $q\rightarrow 1$, $s_{4}(h,k;1)$ is denoted the
classical Hardy-Berndt sums, which is given (\ref{Eq-5}).

By using (\ref{b3-1}), $q$-Hardy-Berndt type sums' $s_{5}(h,k;q)$ is defined
by means of the generating function, $Y_{5}(h,k;q)$:%
\begin{equation}
Y_{5}(h,k;q)=\sum_{%
\begin{array}{c}
m=1 \\ 
2m-1\not\equiv 0\func{mod}k%
\end{array}%
}^{\infty }\frac{F(-\frac{(2m-1)i\pi h}{2k},q)-F(\frac{(2m-1)i\pi h}{2k},q)}{%
2m-1}  \label{s5(h,k)}
\end{equation}%
where $h$ and $k$ are coprime positive integers with $h\geq 1$.

By using (\ref{s5(h,k)}), (\ref{b3-1}) and (\ref{Eq-6}), we obtain the
fallowing theorem.

\begin{theorem}
\label{Theorem 12}Let $h$ and $k$ be denote relatively prime integers with $%
k\geq 1$. If $h$ and $k$ are odd, then%
\begin{equation*}
s_{5}(h,k;q)=\frac{2}{\pi i}Y_{5}(h,k;q).
\end{equation*}
\end{theorem}

\begin{corollary}
Let $h$ and $k$ be denote relatively prime integers with $k\geq 1$. If $h$
and $k$ are odd, then 
\begin{eqnarray*}
s_{5}(h,k;1) &=&\frac{1}{\pi i}Y_{5}(h,k;1) \\
&=&\frac{2}{\pi i}\sum_{%
\begin{array}{c}
m=1 \\ 
2m-1\not\equiv 0\func{mod}k%
\end{array}%
}^{\infty }\frac{1}{m}(\frac{\exp (\frac{\pi ih(2m-1)}{2k})}{1+\exp (\frac{%
\pi ih(2m-1)}{2k})} \\
&&-\frac{\exp (-\frac{\pi ih(2m-1)}{2k})}{1+\exp (-\frac{\pi ih(2m-1)}{2k})}%
).
\end{eqnarray*}
\end{corollary}

Observe that the case $q\rightarrow 1$, $s_{5}(h,k;1)$ is denoted the
classical Hardy-Berndt sums, which is given (\ref{Eq-6}).

\section{$q$-Hardy-Berndt type sums attached to Dirichlet character}

In this section we define new generating functions which are generalizations
of (\ref{b4-1}), (\ref{S1(h,k)}), (\ref{s2(h,k)}), (\ref{s3(h,k)}), (\ref%
{s4(h,k)}) and (\ref{s5(h,k)}). By using (\ref{b3-5}) and%
\begin{equation}
f_{\chi }(t,q)=\sum_{n=1}^{\infty }\chi (n)q^{-n}\exp (-q^{n}[n]t)\text{,
cf. \cite{SimsekJMAA}}  \label{fc(t)}
\end{equation}%
we establish a general theorems related to $q$-Hardy-Berndt type sums with
attached to Dirichlet character.

\begin{align*}
Y_{0,\chi }(h,k;q)& =\sum_{m=1}^{\infty }m^{-p}\sum_{n=1}^{\infty }\chi
(n)(\exp (\frac{q^{-n}[n]\pi ih(2m-1)}{2k}) \\
& -\exp (-\frac{q^{-n}[n]\pi ih(2m-1)}{2k})) \\
& =2i\sum_{m=1}^{\infty }\sum_{n=1}^{\infty }m^{-p}\chi (n)q^{-n}\sin (\frac{%
q^{-n}[n]\pi h(2m-1)}{2k}).
\end{align*}%
Thus we arrive at the following theorem:

\begin{theorem}
\label{Theorem 13}Let $h$ and $k$ be denote relatively prime integers with $%
k\geq 1$. If $h+k$ is odd, then%
\begin{equation*}
Y_{0,\chi }(h,k;q)=2i\sum_{m=1}^{\infty }\sum_{n=1}^{\infty }m^{-p}\chi
(n)q^{-n}\sin (\frac{q^{-n}[n]\pi hm}{2k}).
\end{equation*}
\end{theorem}

Observe that%
\begin{equation*}
\lim_{q\rightarrow 1}Y_{0,\chi }(h,k;q)=2i\sum_{m=1}^{\infty
}\sum_{n=1}^{\infty }m^{-p}\chi (n)\sin (\frac{\pi hmn}{2k}).
\end{equation*}

\begin{proof}[Proof of Theorem 3]
The proof of this theorem is similar to that of Theorem 2. Applying the
generating function to the equation (\ref{A8}), and using (\ref{b3-5}) and (%
\ref{Eq-1}), we arrive at the desired result.
\end{proof}

\begin{remark}
Observe that when $\chi \equiv 1$, principal character and and $q\rightarrow
1$, (\ref{bh0}) reduces to (\ref{A5}) and $S_{1}(h,k;1)$ reduces to $S(h,k)$%
, which is given in (\ref{Eq-1}). If $\chi \equiv 1$, principal character,
then Theorem 3 and Theorem 13 reduce to Theorem 2 and Theorem 7,
respectively.
\end{remark}

We now define generating function of the generalized $q$-Hardy-Bernd type
sums $s_{j,\chi }(h,k;q)$, $j=1,2,3,4,5$. The definition of $q$-analogous of 
$s_{j,\chi }(h,k)$, $j=1,2,3,4,5$ follow precisely along the same lines as
the definition of the $q$-Dedekind type sums, which was established by the
author\cite{SimsekJMAA}. $q$-Hardy-Bernd type sums $s_{1,\chi }(h,k;q)$ is
defined by means of the generating function, $Y_{1,\chi }(h,k;q)$:%
\begin{equation}
Y_{1,\chi }(h,k;q)=\sum_{%
\begin{array}{c}
m=1 \\ 
2m-1\not\equiv 0\func{mod}k%
\end{array}%
}^{\infty }\frac{f_{\chi }(-\frac{(2m-1)i\pi h}{2k},q)-f_{\chi }(\frac{%
(2m-1)i\pi h}{2k},q)}{2m-1},  \label{S11(h,k)}
\end{equation}%
where $f_{\chi }(t,q)$ is defined by (\ref{fc(t)}).

By using (\ref{S11(h,k)}), (\ref{fc(t)}) and (\ref{Eq-2}), we obtain the
fallowing theorem.

\begin{theorem}
\label{Theorem 14}Let $h$ and $k$ be denote relatively prime integers with $%
k\geq 1$. If $h$ is even and $k$ is odd, then%
\begin{equation*}
s_{1,\chi }(h,k;q)=-\frac{2}{\pi i}Y_{1,\chi }(h,k;q).
\end{equation*}
\end{theorem}

\begin{corollary}
Let $h$ and $k$ be denote relatively prime integers with $k\geq 1$. If $h$
is even and $k$ is odd, then%
\begin{eqnarray*}
s_{1,\chi }(h,k;1) &=&-\frac{2}{\pi i}Y_{1,\chi }(h,k,1) \\
&=&-\frac{2}{\pi i}\sum_{%
\begin{array}{c}
m=1 \\ 
2m-1\not\equiv 0\func{mod}k%
\end{array}%
}^{\infty }\sum_{n=1}^{\infty }\frac{\chi (n)}{2m-1}(\frac{\exp (\frac{\pi
ih(2m-1)n}{2k})}{1-\exp (\frac{\pi ih(2m-1)n}{2k})} \\
&&-\frac{\exp (-\frac{\pi ih(2m-1)n}{2k})}{1-\exp (-\frac{\pi ih(2m-1)n}{2k})%
}).
\end{eqnarray*}
\end{corollary}

\begin{remark}
The case $\chi \equiv 1$, principal character and $q\rightarrow 1$, $%
s_{1,1}(h,k;1)$ is denoted the classical Hardy-Berndt sums, which is given (%
\ref{Eq-2}). If $\chi \equiv 1$, principal character, then Theorem 14
reduces to Theorem 8.
\end{remark}

By using (\ref{b3-1}), $q$-Hardy-Berndt type sums' $s_{2,\chi }(h,k;q)$ is
defined by means of the generating function, $Y_{2,\chi }(h,k;q)$:%
\begin{equation}
Y_{2,\chi }(h,k;q)=\sum_{%
\begin{array}{c}
m=1 \\ 
2m\not\equiv 0\func{mod}k%
\end{array}%
}^{\infty }\frac{F_{\chi }(-\frac{mi\pi h}{k},q)-F_{\chi }(\frac{mi\pi h}{k}%
,q)}{m}  \label{s22(h,k)}
\end{equation}%
By using (\ref{s22(h,k)}), (\ref{b3-5}) and (\ref{Eq-3}), we obtain the
fallowing theorem.

\begin{theorem}
\label{Theorem 15}Let $h$ and $k$ be denote relatively prime integers with $%
k\geq 1$. If $h$ is odd and $k$ is even, then%
\begin{equation*}
s_{2,\chi }(h,k;q)=-\frac{1}{2\pi i}Y_{2,\chi }(h,k;q).
\end{equation*}
\end{theorem}

\begin{corollary}
Let $h$ and $k$ be denote relatively prime integers with $k\geq 1$. If $h$
is odd and $k$ is even, then%
\begin{eqnarray*}
s_{2,\chi }(h,k;1) &=&-\frac{1}{2\pi i}Y_{2,\chi }(h,k,1) \\
&=&-\frac{2}{\pi i}\sum_{%
\begin{array}{c}
m=1 \\ 
2m\not\equiv 0\func{mod}k%
\end{array}%
}^{\infty }\sum_{n=1}^{\infty }\frac{\chi (n)}{m}(\frac{\exp (\frac{\pi ihmn%
}{k})}{1+\exp (\frac{\pi ihmn}{k})} \\
&&-\frac{\exp (-\frac{\pi ihmn}{k})}{1+\exp (-\frac{\pi ihmn}{k})}).
\end{eqnarray*}
\end{corollary}

\begin{remark}
The case$\chi \equiv 1$, principal character and $q\rightarrow 1$, $%
s_{2,1}(h,k;1)$ is denoted the classical Hardy-Berndt sums, which is given (%
\ref{Eq-3}). If $\chi \equiv 1$, principal character, then Theorem 15
reduces to Theorem 9.
\end{remark}

By using (\ref{b3-5}), $q$-Hardy-Berndt type sums' $s_{3,\chi }(h,k;q)$ is
defined by means of the generating function, $Y_{3,\chi }(h,k;q)$:%
\begin{equation}
Y_{3,\chi }(h,k;q)=\sum_{m=1}^{\infty }\frac{F_{\chi }(-\frac{mi\pi h}{k}%
,q)-F_{\chi }(\frac{mi\pi h}{k},q)}{m}  \label{s33(h,k)}
\end{equation}%
By using (\ref{s33(h,k)}), (\ref{b3-5}) and (\ref{Eq-4}), we obtain the
fallowing theorem.

\begin{theorem}
\label{Theorem 16}Let $h$ and $k$ be denote relatively prime integers with $%
k\geq 1$. If $k$ is odd, then%
\begin{equation*}
s_{3,\chi }(h,k;q)=\frac{1}{\pi i}Y_{3,\chi }(h,k;q).
\end{equation*}
\end{theorem}

\begin{corollary}
Let $h$ and $k$ be denote relatively prime integers with $k\geq 1$. If $k$
is odd, then 
\begin{eqnarray*}
s_{3,\chi }(h,k;1) &=&\frac{1}{\pi i}Y_{3,\chi }(h,k;1) \\
&=&\frac{1}{\pi i}\sum_{m=1}^{\infty }\sum_{n=1}^{\infty }\frac{\chi (n)}{m}(%
\frac{\exp (\frac{\pi ihmn}{k})}{1+\exp (\frac{\pi ihmn}{k})} \\
&&-\frac{\exp (-\frac{\pi ihmn}{k})}{1+\exp (-\frac{\pi ihmn}{k})}).
\end{eqnarray*}
\end{corollary}

\begin{remark}
The case $\chi \equiv 1$, principal character and $q\rightarrow 1$, $%
s_{3,1}(h,k;1)$ is denoted the classical Hardy-Berndt sums, which is given (%
\ref{Eq-4}). If $\chi \equiv 1$, principal character, then Theorem 16
reduces to Theorem 10.
\end{remark}

$q$-Hardy-Bernd type sums $s_{4,\chi }(h,k;q)$ is defined by means of the
generating function, $Y_{4,\chi }(h,k,q)$:%
\begin{equation}
Y_{4,\chi }(h,k;q)=\sum_{m=1}^{\infty }\frac{f_{\chi }(-\frac{(2m-1)i\pi h}{%
2k},q)-f_{\chi }(\frac{(2m-1)i\pi h}{2k},q)}{2m-1},  \label{s44(h,k)}
\end{equation}%
where $f_{\chi }(t,q)$ is defined by (\ref{fc(t)}).

By using (\ref{s44(h,k)}), (\ref{fc(t)}) and (\ref{Eq-5}), we obtain the
fallowing theorem.

\begin{theorem}
\label{Theorem 17}Let $h$ and $k$ be denote relatively prime integers with $%
k\geq 1$. If $h$ is odd, then%
\begin{equation*}
s_{4,\chi }(h,k;q)=\frac{4}{\pi i}Y_{4,\chi }(h,k;q).
\end{equation*}
\end{theorem}

\begin{corollary}
Let $h$ and $k$ be denote relatively prime integers with $k\geq 1$. If $h$
is odd, then%
\begin{eqnarray*}
s_{4,\chi }(h,k;1) &=&\frac{4}{\pi i}Y_{4,\chi }(h,k;1) \\
&=&\frac{4}{\pi i}\sum_{m=1}^{\infty }\sum_{n=1}^{\infty }\frac{\chi (n)}{%
2m-1}(\frac{\exp (\frac{\pi ihn(2m-1)}{2k})}{1-\exp (\frac{\pi ihn(2m-1)}{2k}%
)} \\
&&-\frac{\exp (-\frac{\pi ihn(2m-1)}{2k})}{1-\exp (-\frac{\pi ihn(2m-1)}{2k})%
}).
\end{eqnarray*}
\end{corollary}

\begin{remark}
The case $\chi \equiv 1$, principal character and $q\rightarrow 1$, $%
s_{4,1}(h,k;1)$ is denoted the classical Hardy-Berndt sums, which is given (%
\ref{Eq-5}). If $\chi \equiv 1$, principal character, then Theorem 17
reduces to Theorem 11.
\end{remark}

By using (\ref{b3-5}), $q$-Hardy-Berndt type sums' $s_{5,\chi }(h,k;q)$ is
defined by means of the generating function, $Y_{5,\chi }(h,k;q)$:%
\begin{equation}
Y_{5,\chi }(h,k;q)=\sum_{%
\begin{array}{c}
m=1 \\ 
2m-1\not\equiv 0\func{mod}k%
\end{array}%
}^{\infty }\frac{F_{\chi }(-\frac{(2m-1)i\pi h}{2k},q)-F_{\chi }(\frac{%
(2m-1)i\pi h}{2k},q)}{2m-1}  \label{s55(h,k)}
\end{equation}%
where $h$ and $k$ are coprime positive integers with $h\geq 1$.

By using (\ref{s55(h,k)}), (\ref{b3-5}) and (\ref{Eq-6}), we obtain the
fallowing theorem.

\begin{theorem}
\label{Theorem 18}Let $h$ and $k$ be denote relatively prime integers with $%
k\geq 1$. If $h$ and $k$ are odd, then%
\begin{equation*}
s_{5,\chi }(h,k;q)=\frac{2}{\pi i}Y_{5,\chi }(h,k,q).
\end{equation*}
\end{theorem}

\begin{corollary}
Let $h$ and $k$ be denote relatively prime integers with $k\geq 1$. If $h$
and $k$ are odd, then 
\begin{eqnarray*}
s_{5,\chi }(h,k;1) &=&\frac{1}{\pi i}Y_{5,\chi }(h,k,1) \\
&=&\frac{2}{\pi i}\sum_{%
\begin{array}{c}
m=1 \\ 
2m-1\not\equiv 0\func{mod}k%
\end{array}%
}^{\infty }\sum_{n=1}^{\infty }\frac{\chi (n)}{2m-1}(\frac{\exp (\frac{\pi
ihn(2m-1)}{2k})}{1+\exp (\frac{\pi ihn(2m-1)}{2k})} \\
&&-\frac{\exp (-\frac{\pi ihn(2m-1)}{2k})}{1+\exp (-\frac{\pi ihn(2m-1)}{2k})%
}).
\end{eqnarray*}
\end{corollary}

\begin{remark}
The case $\chi \equiv 1$, principal character and $q\rightarrow 1$, $%
s_{5,1}(h,k;1)$ is denoted the classical Hardy-Berndt sums, which is given (%
\ref{Eq-6}). If $\chi \equiv 1$, principal character, then Theorem 18
reduces to Theorem 12.
\end{remark}

\section{New generating functions related to the Riemann zeta function, $L$%
-function and Genocchi zeta function}

In this section, our goal is to define new generating functions. By applying
the Mellin transformation to these functions, we can give some new relations
which are related to Riemann zeta functions, $q$-Riemann zeta function, $q$-$%
L$-function and $q$-Genocchi zeta function.

We define%
\begin{equation}
y_{0}(t,q)=\sum_{m=1}^{\infty }\frac{F(-(2m-1)it,q)-F((2m-1)it,q)}{2m-1},
\label{yS0}
\end{equation}%
where $F(t,q)$ is defined by the relation (\ref{b3-1}).

\begin{theorem}
\label{Theorem 19}We have%
\begin{equation*}
\frac{1}{\Gamma (s)}\int_{0}^{\infty
}t^{s-1}y_{0}(t,q)dt=(i)^{-s}((-1)^{-s}-1)\Im _{G,q}(s)\zeta ^{\ast }(s+1),
\end{equation*}%
where $\Im _{G,q}(s)$ is the \textit{Genocchi} zeta functions and%
\begin{equation*}
\zeta ^{\ast }(s):=\sum_{m=1}^{\infty }\frac{1}{(2m-1)^{s}}.
\end{equation*}
\end{theorem}

\begin{proof}
By applying the Mellin transformation to the equation (\ref{yS0}) and using (%
\ref{b3-1}), we find that%
\begin{align*}
& \frac{1}{\Gamma (s)}\int_{0}^{\infty }t^{s-1}y_{0}(t,q)dt \\
& =\frac{1}{\Gamma (s)}\int_{0}^{\infty }t^{s-1}\left( \sum_{m=1}^{\infty }%
\frac{F(-(2m-1)it,q)-F((2m-1)it,q)}{2m-1}\right) dt \\
& =\frac{1}{\Gamma (s)}\sum_{m=1}^{\infty }\frac{1}{2m-1}(\sum_{n=1}^{\infty
}(-1)^{n}q^{-n}\int_{0}^{\infty }t^{s-1}\exp \left( q^{-n}[n](2m-1)ti\right)
dt \\
& -\sum_{n=1}^{\infty }(-1)^{n}q^{-n}\int_{0}^{\infty }t^{s-1}\exp \left(
-q^{-n}[n](2m-1)ti\right) dt) \\
& =\frac{1}{\Gamma (s)}\sum_{m=1}^{\infty }\frac{(i)^{-s}((-1)^{-s}-1)}{%
(2m-1)^{s+1}}\left( \sum_{n=1}^{\infty }\frac{(-1)^{n}q^{-n}}{(q^{-n}[n])^{s}%
}\int_{0}^{\infty }u^{s-1}e^{-u}du\right) \\
& =(i)^{-s}((-1)^{-s}-1)\sum_{m=1}^{\infty }\frac{1}{(2m-1)^{s+1}}%
\sum_{n=1}^{\infty }\frac{(-1)^{n}q^{-n}}{(q^{-n}[n])^{s}}.
\end{align*}%
After some elementary calculations, we obtain the desired result.
\end{proof}

\begin{remark}
In \cite{srivastava and choi}, p. 96, Eq. 2.3 (1), Srivastava and Choi
defined Riemann zeta function as follows:%
\begin{equation*}
\zeta (s)=\left\{ 
\begin{array}{c}
\sum_{m=1}^{\infty }n^{-s}=\frac{1}{1-2^{-s}}\sum_{m=1}^{\infty }(2m-1)^{-s},%
\text{ \ \ }\func{Re}(s)>1 \\ 
(1-2^{1-s})^{-1}\sum_{m=1}^{\infty }(-1)^{m-1}m^{-s},\text{ \ \ (}\func{Re}%
(s)>0\text{, }s\neq 1.%
\end{array}%
\right.
\end{equation*}%
For $\func{Re}(s)>1,$ by using this definition, relation between $\zeta
^{\ast }(s)$ and $\zeta (s)$ is given by 
\begin{equation*}
\zeta ^{\ast }(s)=(1-2^{-s})\zeta (s).
\end{equation*}%
In \cite{srivastava and choi}, p. 121, Eq. 2.5 (1), They defined
Hurwitz-Lerch zeta functions as follows%
\begin{equation*}
\Phi (z,s,a):=\sum_{m=1}^{\infty }\frac{z^{m}}{(m+a)^{s}},
\end{equation*}%
where $a\in \mathbb{C}\backslash \mathbb{Z}_{0}^{-}$; $s\in \mathbb{C}$ when 
$\mid z\mid <1$; $\func{Re}(s)>1$ when $\mid z\mid =1$. In \cite{srivastava
and choi}, p. 339, Eq. 6.1 (18), They gave the following identity%
\begin{equation*}
\sum_{m=1}^{\infty }\frac{z^{m}}{(2m-1)^{s}}=(2b)^{-s}\sum_{j=1}^{b}\Phi
(z^{b},s,\frac{2j-1}{2b})z^{j-1}.
\end{equation*}%
By substituting $\func{Re}(s)>1$ when $z=1$ into the above identity, we have%
\begin{eqnarray*}
\zeta ^{\ast }(s) &=&(2b)^{-s}\sum_{j=1}^{b}\Phi (1,s,\frac{2j-1}{2b}) \\
&=&(2b)^{-s}\zeta (s,\frac{2j-1}{2b}),
\end{eqnarray*}%
where $\zeta (s,x)$ denotes Hurwitz zeta functions.
\end{remark}

We define the following generating function:%
\begin{equation}
y_{1}(t,q)=\sum_{m=1}^{\infty }\frac{f(-(2m-1)it,q)-f((2m-1)it,q)}{2m-1},
\label{ys1}
\end{equation}%
where $f(t,q)$ is defined by (\ref{f(t,q)}).

\begin{theorem}
\label{Theorem 20}We have%
\begin{equation*}
\frac{1}{\Gamma (s)}\int_{0}^{\infty
}t^{s-1}y_{1}(t,q)dt=(i)^{-s}((-1)^{-s}-1)\zeta _{q}(s)\zeta ^{\ast }(s+1),
\end{equation*}%
where $\zeta _{q}(s)$ is the $q$-Riemann zeta function, which is defined in 
\cite{SimsekJMAA}.
\end{theorem}

\begin{proof}
By applying the Mellin transformation to the equation (\ref{ys1}) and using (%
\ref{f(t,q)}), we find that%
\begin{eqnarray*}
&&\frac{1}{\Gamma (s)}\int_{0}^{\infty }t^{s-1}y_{1}(t,q)dt \\
&=&\frac{1}{\Gamma (s)}\int_{0}^{\infty }t^{s-1}\left( \sum_{m=1}^{\infty }%
\frac{f(-(2m-1)it,q)-f((2m-1)it,q)}{2m-1}\right) dt.
\end{eqnarray*}%
The remainder of the proof runs parallel to that of \ Theorem 19, so we
choose to omit the details involved.
\end{proof}

We define the following generating function:

\begin{equation*}
y_{2}(t,q)=\sum_{m=1}^{\infty }\frac{F(-tmi,q)-F(tmi,q)}{m}.
\end{equation*}%
By applying the Mellin transformation to the above function, we arrive at
the following theorem:

\begin{theorem}
\label{Theorem 21}We have%
\begin{equation*}
\frac{1}{\Gamma (s)}\int_{0}^{\infty
}t^{s-1}y_{2}(t,q)dt=(i)^{-s}((-1)^{-s}-1)\Im _{G,q}(s)\zeta (s+1),
\end{equation*}%
where $\zeta (s+1)$ is the Riemann zeta functions.
\end{theorem}

\begin{remark}
Proof of Theorem 21 run parallel to that of Theorem 19, so we choose to omit
the details involved.
\end{remark}

We define the following generating function attached to Drichlet character:%
\begin{equation}
y_{0,\chi }(t,q)=\sum_{m=1}^{\infty }\frac{F_{\chi }(-(2m-1)it,q)-F_{\chi
}((2m-1)it,q)}{2m-1},  \label{ys0x}
\end{equation}%
where $F_{\chi }(t,q)$ is defined by (\ref{fc(t)}).

\begin{theorem}
\label{Theorem 22}We have%
\begin{equation*}
\frac{1}{\Gamma (s)}\int_{0}^{\infty }t^{s-1}y_{0,\chi
}(t,q)dt=(i)^{-s}((-1)^{-s}-1)l_{G,q}(s,\chi )\zeta ^{\ast }(s+1).
\end{equation*}
\end{theorem}

\begin{proof}
By applying the Mellin transformation to the equation (\ref{ys0x}) and using
(\ref{fc(t)}), we find that%
\begin{align*}
& \frac{1}{\Gamma (s)}\int_{0}^{\infty }t^{s-1}y_{0,\chi }(t,q)dt \\
& =\frac{1}{\Gamma (s)}\int_{0}^{\infty }t^{s-1}\left( \sum_{m=1}^{\infty }%
\frac{F_{\chi }(-(2m-1)it,q)-F_{\chi }((2m-1)it,q)}{2m-1}\right) dt \\
& =\frac{1}{\Gamma (s)}\sum_{m=1}^{\infty }\frac{1}{2m-1}(\sum_{n=1}^{\infty
}(-1)^{n}\chi (n)q^{-n}\int_{0}^{\infty }t^{s-1}\exp \left(
q^{-n}[n](2m-1)ti\right) dt \\
& -\sum_{n=1}^{\infty }(-1)^{n}\chi (n)q^{-n}\int_{0}^{\infty }t^{s-1}\exp
\left( -q^{-n}[n](2m-1)ti\right) dt) \\
& =\frac{1}{\Gamma (s)}\sum_{m=1}^{\infty }\frac{(i)^{-s}((-1)^{-s}-1)}{%
(2m-1)^{s+1}}\left( \sum_{n=1}^{\infty }\frac{(-1)^{n}\chi (n)q^{-n}}{%
(q^{-n}[n])^{s}}\int_{0}^{\infty }u^{s-1}e^{-u}du\right) \\
& =(i)^{-s}((-1)^{-s}-1)\sum_{m=1}^{\infty }\frac{1}{(2m-1)^{s+1}}%
\sum_{n=1}^{\infty }\frac{(-1)^{n}\chi (n)q^{-n}}{(q^{-n}[n])^{s}}.
\end{align*}%
After some elementary calculations, we obtain the desired result.
\end{proof}

We define the following generating function:%
\begin{equation}
y_{1,\chi }(t,q)=\sum_{m=1}^{\infty }\frac{f_{\chi }(-(2m-1)it,q)-f_{\chi
}((2m-1)it,q)}{2m-1},  \label{ys1x}
\end{equation}%
where $f_{\chi }(t,q)$ is defined by (\ref{fc(t)}).

\begin{theorem}
\label{t23}We have%
\begin{equation*}
\frac{1}{\Gamma (s)}\int_{0}^{\infty }t^{s-1}y_{1,\chi
}(t,q)dt=(i)^{-s}((-1)^{-s}-1)L_{q}(s,\chi )\zeta ^{\ast }(s+1),
\end{equation*}%
where $L_{q}(s,\chi )$ is the $\mathit{q}$-$L$ function, which is defined in 
\cite{SimsekJMAA}.
\end{theorem}

\begin{proof}
By applying the Mellin transformation to the equation (\ref{ys1x}) and using
(\ref{fc(t)}), we find that%
\begin{eqnarray*}
&&\frac{1}{\Gamma (s)}\int_{0}^{\infty }t^{s-1}y_{1,\chi }(t,q)dt \\
&=&\frac{1}{\Gamma (s)}\int_{0}^{\infty }t^{s-1}\left( \sum_{m=1}^{\infty }%
\frac{f_{\chi }(-(2m-1)it,q)-f_{\chi }((2m-1)it,q)}{2m-1}\right) dt.
\end{eqnarray*}%
The remainder of the proof runs parallel to that of Theorem 22, so we choose
to omit the details involved.
\end{proof}

\begin{remark}
By setting $t=\frac{\pi ih}{k}$, $h$, $k\in \mathbb{Z}$ with $(h,k)=1$ in (%
\ref{yS0}), (\ref{ys1}), (\ref{ys0x}) and (\ref{ys1x}), we obtain generating
functions in Section 4 and 5, respectively.
\end{remark}

\begin{acknowledgement}
This paper was supported by the Scientific Research Project Administration
Akdeniz University.
\end{acknowledgement}

\end{document}